\documentclass[10pt,a4paper,article]{amsart}
\usepackage{amssymb,stmaryrd}
\usepackage{amsthm,amstext}
\usepackage{amssymb,amsmath}
\usepackage{amsfonts}
\usepackage{MnSymbol,slashed}
\usepackage{fancyhdr}
\usepackage{graphicx}
\usepackage{pdfpages}
\usepackage{float}
\usepackage{subfig}
\usepackage{simplewick}
\usepackage{graphics,psfrag}
\usepackage{pstricks}
\usepackage{latexsym}
\usepackage{changepage}
\usepackage{caption}
\usepackage{paralist}

\theoremstyle{plain}
\newtheorem{theorem}{Theorem}[section] 

\newtheorem{lemma}[theorem]{Lemma}

\newtheorem{proposition}[theorem]{Proposition}

\newcommand{\fun}[3]{#1: #2 \rightarrow #3}
\newcommand{\onef}[1]{\Omega^{#1}}
\newcommand{\tchrs}[2]{{\Gamma}^{#1}_{\hphantom{#1}#2}}

\newcommand{\dirac}{{ \slashed{D} }}
\newcommand*{\doublenabla}{%
  \nabla\mkern-12mu\nabla}

\newcommand{\CH}{\hbox{{$\mathcal H$}}}
\newcommand{\CS}{{\hbox{{$\mathcal S$}}}}

\newcommand{\CJ}{\hbox{{$\mathcal J$}}}

\newcommand{\C}{\mathbb{C}}
\newcommand{\R}{\mathbb{R}}

\newcommand{\Z}{\mathbb{Z}}
\newcommand{\N}{\mathbb{N}}

\newcommand{\del}{\partial}

\newcommand{\extd}{\mathrm{d}}

\newcommand{\eps}{{\epsilon}}
\newcommand{\tens}{\mathop{{\otimes}}}
\newcommand{\la}{{\triangleright}}

\newcommand{\id}{\mathrm{id}}

\newcommand{\<}{\langle}
\renewcommand{\>}{\rangle}

\begin{document}

\author{Evelyn Lira-Torres and Shahn Majid}
\address{School of Mathematical Sciences\\ Queen Mary University of London \\ Mile End Rd, London E1 4NS }
\email{ e.y.liratorres@qmul.ac.uk, s.majid@qmul.ac.uk}
\thanks{Ver. 2.2 The first author was partially supported by CONACyT and Fundaci\'on Alberto y Dolores Andrade (M\'exico)}

\title{Geometric Dirac operator on the fuzzy sphere}
	\begin{abstract} We construct a Connes spectral triple or `Dirac operator' on the non-reduced fuzzy sphere $\C_\lambda[S^2]$ as realised using quantum Riemannian geometry with a central quantum metric $g$ of Euclidean signature and its associated quantum Levi-Civita connection. The Dirac operator is characterised uniquely up to unitary equivalence within our quantum Riemannian geometric setting and an assumption that the spinor bundle is trivial and rank 2 with a central basis. The spectral triple has KO dimension 3 and in the case of the round metric, essentially recovers a previous proposal motivated by rotational symmetry. 
	\end{abstract}
\keywords{Spectral triple, fuzzy sphere, quantum geometry, fuzzy monopole, noncommutative geometry, angular momentum algebra, coadjoint quantisation}
\maketitle

\section{Introduction}

Noncommutative geometry, or the idea that coordinate algebras can be noncommutative, can take many forms but one which has been influential is a `top down' approach in which all of the geometry is encoded in Connes' notion of a {\em spectral triple}\cite{Con,ConMar}. This is an axiomatic framework for a `Dirac operator' $\dirac$ for a possibly noncommutative $*$-algebra $A$ acting on a Hilbert space $\CH$ where $A$ is also represented. Meanwhile, in the last decades, there has also emerged a complementary `bottom up' approach to quantum Riemannian geometry, e.g. \cite{BegMa2, BegMa, ArgMa1, ArgMa, LirMa, MaPac} in which we start with a choice of differential bimodule $(\Omega^1,\extd)$ on $A$, typically guided by classical or quantum group symmetry, and then proceed to build up the different layers of geometry starting with the quantum metric $g\in\Omega^1\tens_A\Omega^1$ and quantum Levi-Civita connection $\nabla:\Omega^1\to \Omega^1\tens_A\Omega^1$, subject to various axioms. Principal among these in the strongest version is existence of a bimodule inverse metric $(\ ,\ ):\Omega^1\tens_A\Omega^1\to A$ and a requirement  that $\nabla$ is a bimodule connection in the sense of\cite{DVM,Mou}. There are also weaker versions and variants, but it has made sense to explore the theory in this nicest case first. This more constructive approach may not necessarily reach the level of  Dirac operator at all and  grew out of experience with quantum group-symmetric models where the axioms of a spectral triple  were not necessarily a good fit. It was shown in \cite{BegMa:spe}, however, that this it is sometimes possible to arrive at something close to a spectral triple  in this way. This possibility of a conjunction of the `top down' and `bottom up' approaches raises the interesting problem in general of identifying which Connes spectral triples on a given algebra can be quantum geometrically realised in this way. Such models would then have both the deep significance of Connes' approach from the point of view of KO homology, but also explicit access to the various layers of quantum geometry that could be useful in mathematical physics. Requiring  this could, for example, help single out particular Connes spectral triples and this could then have predictive implications, for example in the Standard Model where a finite spectral triple tensored onto spacetime can be used to encode the structure of elementary particles\cite{Con0,ConMar,DabSit}. It could also single out examples of interest for the modelling of quantum gravity effects, a general idea that has a long history e.g. \cite{Sny,Ma:pla} and which is currently visible in 2+1 quantum gravity e.g. \cite{Hoo, FreMa, MaSch}. 

In this paper, we explore this conjunction in the case of `fuzzy spheres'. These have been extensively studied since they were introduced in \cite{Mad}, for example \cite{Gro,Gro1, Gro2, Bar, Wat,BalPad,And}. Some of these works construct spectral triples of some kind, but as matrix algebras or `finite noncommutative geometries'. Here we work on the  non-reduced fuzzy sphere $A=\C_\lambda[S^2]$, which is an infinite-dimensional algebra obtained as the coadjoint quantisation of a standard sphere, i.e. the angular momentum algebra $U(su_2)$ with fixed value of the quadratic Casimir so as to deform a unit sphere. It's quantum Riemannian geometry for a central quantum metric $g$ was recently studied in \cite{LirMa} and in this sequel we will supply a canonical Dirac operator and spectral triple. Specifically, $\C_\lambda[S^2]$ has generators $x_i$ obeying
\begin{equation}\label{fuzzyalg} [x_i,x_j]=2\imath\lambda_p \eps_{ijk}x_k,\quad  \sum_i x_i^2=1-\lambda_p^2,\end{equation}
for $i,j,k\in{1,2,3}$, where $\lambda_p$ is a deformation parameter. We sum over repeated indices and $\eps_{ijk}$ is the totally antisymmetric tensor.  If we scale our generators by a factor $L$ so as to deform a sphere or radius $L$ and if such models arise from quantum gravity, then one might expect that $\lambda_p\sim l_p/L$ where $l_p$ is the Planck scale. The result of our analysis is that we are led, see Proposition~\ref{propdirg}, to
\begin{equation}\label{diracg} (\dirac\psi)_\alpha =\del_i\psi_\alpha C^{j\alpha}{}_\beta+{\imath\over 4}{{\rm Tr}(g)\over\sqrt{\det(g)}}\psi_\alpha\end{equation}
where $C^i$ are Clifford algebra matrices for the metric tensor $g$ in a basis of a certain 3D differential calculus $\Omega^1$, $\del_i$ are defined by the exterior derivative in the basis and $\psi_\alpha$ is a 2-component spinor. The second term is a curvature in the same spirit as the Ricci curvature in \cite{LirMa}, which also involved determinants and traces of $g$ (but a trace of $g^2$ not just of $g$ as now). We assume for simplicity that $g$ has Euclidean signature, of which the simplest case is the `round metric' $g_{ij}=\delta_{ij}$, in which case we can take $C^i=\sigma^i$ the Pauli matrices. 

In fact, the  relations (\ref{fuzzyalg})  are arranged so that if $\lambda_p=1/n$ for $n\in \N$ then the $n$-dimensional representation of $U(su_2)$ descends to $A$. It has a large kernel but if we were to quotient out by this kernel then we would have something isomorphic to $M_n(\C)$. These `reduced' fuzzy spheres are the ones studied in preceding works as examples of finite noncommutative geometries\cite{Mad,Gro,Gro2,Wat,Bar,BalPad,And}. By contrast, $\C_\lambda[S^2]$ deforms the usual functions on a sphere for all $\lambda_p$ including generic values. This, in turn, means that the Hilbert space of spinors on which our $\dirac$ is represented is infinite dimensional, albeit we will be able to work with finite-dimensional subspaces $\CS_l$ of orbital angular momentum $l$. In the case of the round metric $g_{ij}=\delta_{ij}$, the eigenvalues of $\imath\dirac$ are 
\[  \lambda_{l,\pm}=-{1\over 4}\pm(l+{1\over 2}),\quad l=0,1,2,\cdots\]
with only $\lambda_{0,-}$ when $l=0$. The total set of eigenvalues maps by a constant shift of $-1/4$ to the total set of eigenvalues $\pm(l+1)$ of a classical sphere\cite{Cam}, so the main effect is a kind of `zero point' shift. Everything descends to the reduced case when $\lambda_p=1/n$ and in the round metric case we then essentially recover a proposal of D'Andrea, Lizzi and Varilly\cite{And} coming very differently from rotational symmetry and Berezin quantisation considerations and without the $-1/4$ shift. This in turn was built on a Dirac operator first proposed in \cite{BalPad} motivated by rotational symmetry. In our case, $\dirac$ comes uniquely up to equivalence out of the quantum geometry and an assumption of a trivial spinor bundle as an $A$-bimodule with central basis, and is rotationally invariant just because the round metric $g$ is. This uniqueness is in Proposition~\ref{dirU}. 

Section~\ref{secpre} has some very minimal preliminaries on the formalism of quantum Riemannian geometry on the one hand and of Connes' formalism on the other, to fix notation. We refer to \cite{BegMa,Con,ConMar} for more details of the formalisms. Section~\ref{secsetup} then sets up the equations to be solved in the case of a trivial spinor bundle with central basis. Section~\ref{secround} solves them to construct $\dirac$.  We conclude in Section~\ref{secrem} with some  directions for further work. 

\section{Preliminaries: Fuzzy sphere and outline of the formalism} \label{secpre}

Here we outline the steps needed to get to a Dirac operator coming from the bottom up as laid out in \cite{BegMa:spe}. Because this is a sequel to \cite{LirMa}, and given textbooks \cite{BegMa,Con,ConMar}, we only give bare details for orientation and so as to fix notations. We let $A$ be a unital $*$-algebra over $\C$. 
From the constructive `bottom up' point of view, we fix the following structures in order. 

\subsection{}\label{sec21} $(\Omega(A),\extd)$ a differential graded `exterior algebra' over $A$. We require this to be generated by $A,\extd A$ as it would be classically. For Riemannian geometry we only really need $\Omega^1,\Omega^2$ as $A$-$A$ bimodules in the lowest degree of the complex of differential forms. We denote by $\wedge$ the wedge product of forms of degree $\ge 1$. We also ask for $\ker \extd: A\to \Omega^1$ to be spanned by 1, i.e. we are interested in connected differential calculi.  We require $*$ to extend as a graded-involution commuting with $\extd$. Specifically, for $\C_\lambda[S^2]$, we use the calculus\cite[Ex.~1.46]{BegMa}
\[ [s^i , x_j ]=0,\quad \extd x_i= \epsilon_{ijk}x_{j}s^k,\]
where $s^i$, $i=1,2,3$ are  central basis of $\Omega^1$ over $A$. This  projects when $\lambda_p=1/n$ to  similar differentials on the reduced $M_n(\C)$ fuzzy sphere \cite{Mad}. It differs from a classical sphere in having one dimension too many, which one can think of as a `normal direction' $\theta'=x^i\extd x^i/(2\imath\lambda_p)$ forced on us by the requirement of rotational invariance. Here
\[   s^i= \tfrac{1}{(1-\lambda_p^2)}\left( x_i\theta'+\epsilon_{ijk}(\extd x_j) x_k\right)\]
The exterior algebra is given by the $s^i$ mutually anticommuting (a Grassmann algebra) and $\extd s^i=-{1\over 2}\eps_{ijk}s^j\wedge s^k$. The $*$-structure is $x_i^*=x_i$ and $s^i{}^*=s^i$. In the classical limit, one can consider that $\theta'\to 0$ as well as $\lambda_p\to 0$, in which case $s^i$ become geometrically the 1-form versions of the Killing vectors on $S^2$ associated to the rotational symmetry, i.e. to generators of orbital angular momentum. Indeed, the partial derivatives $\del_i$ defined by $\extd f=(\del_i f)s^i$ can be given by $\del_i=[x_i,\ ]/(2\imath\lambda_p)$ in view of the above commutation relations as one might expect for orbital angular momentum in our context. 

\subsection{}\label{sec22} We define a quantum metric as  $g\in \onef{1}\otimes_A \onef{1}$ such that there exists an inverse $(\ ,\ ):\Omega^1\tens_A\Omega^1\to A$ which is a bimodule map. Inverse here means in the usual sense but turns out to require that $g$ is central. One can (and we will) optionally require $g$ to be quantum symmetric in the sense $\wedge(g)=0$.  The proof that $g$ has to be central is in \cite{BegMa2} and \cite[Lemma~1.16]{BegMa}.  We also require ${\rm flip}(*\tens *)(g)=g$ for compatibility with $*$. For the fuzzy sphere, this leads to 
\begin{equation}\label{g} g=g_{ij}s^i\tens s^j,\end{equation}
where $g_{ij}$ is a real symmetric and invertible matrix. This is much more restrictive than on a classical sphere and could be viewed as an important subclass within a more general theory, but sufficient for the `round metric' $g_{ij}=\delta_{ij}$. 

\subsection{}\label{sec23}  We require a bimodule connection $\nabla:\Omega^1\to \Omega^1\tens\Omega^1$ obeying  the left Leibniz rule
\[ \nabla(a.\omega)=\extd a\tens\omega+a.\nabla\omega\]
for all $a\in A$ and $\omega\in \Omega^1$ as usual in noncommutative geometry, and in addition another Leibniz rule 
\[ \nabla (\omega a)=(\nabla \omega)a+\sigma (\omega \otimes\extd  a),\] 
for a bimodule map $\fun{\sigma}{\onef{1}\otimes _A \onef{1}}{\onef{1}\otimes _A \onef{1}}$,  the \textit{generalised braiding}. Such `bimodule connections' on $A$-bimodules\cite{DVM,Mou} extend automatically to tensor products, with the result, in our case, that $\nabla g=0$ makes sense. We say that $\nabla$ is a quantum Levi-Civita connection (QLC) if this holds and if the torsion $T_\nabla=\wedge\nabla-\extd:\Omega^1\to \Omega^2$ vanishes. We also require $\sigma\circ {\rm flip}(*\tens *)\nabla= \nabla \circ *$ for compatibility with $*$. These are a well-studied set of axioms\cite[Chap~8]{BegMa} for which many interesting examples are known, e.g.  \cite{BegMa2, ArgMa1, ArgMa, LirMa, MaPac}. For $\C_\lambda[S^2]$, there is a unique QLC under the reasonable assumption (given the metric above) that its coefficients are constant in the $s^i$ basis, namely\cite{LirMa}
\begin{equation}\label{QLC} \nabla s^i= -\tfrac{1}{2}\tchrs{i}{jk}s^j \otimes s^k,\quad \Gamma_{ijk}=2\eps_{ikm}g_{mj}+{\rm Tr}(g)\eps_{ijk},\end{equation}
where $\Gamma_{ijk}:=g_{im}\Gamma^m{}_{jk}$. For the round metric, this reduces to $\Gamma_{ijk}=\eps_{ijk}$ in keeping with formulae in \cite{Mad} for the reduced $M_n(\C)$ fuzzy sphere.

 \subsection{}\label{sec24} We now progress to describe the extension of the quantum Riemannian geometry  formalism to a Dirac operator according to the steps laid out in \cite{BegMa:spe}. For this, we require a bimodule map `Clifford action' $\la:\Omega^1\tens_A \CS\to \CS$ and a compatible left bimodule connection $\nabla_{\CS}:\CS\to \Omega^1\tens_A\CS$. We then define $\dirac=\la \circ \nabla_{\CS}: \CS\to \CS$ as the associated Dirac operator. 
 
To approach Connes' axioms, we further require an antilinear map $\CJ:\CS\to \CS$ such that \cite{BegMa:spe,BegMa}
\begin{equation}\label{J} (\CJ s).a=\CJ (a^* s),\quad \CJ^2=\eps\, \id\end{equation}
 for all $a\in A$ and $s\in\CS$ and some fixed sign $\eps=\pm 1$. Next, \cite{BegMa:spe,BegMa} requires, in concrete terms, the compatibility
\begin{equation}\label{JnablaS}  \sigma_{\CS} \circ {\rm flip}\circ (*\tens \CJ)\circ \nabla_{\CS} =\nabla_{\CS}\circ\CJ ,\quad  \CJ(\omega\la s)=\eps'\, \la\circ\sigma_{\CS}(\CJ s\tens \omega^*)\end{equation}
for all $s\in \CS$ and $\omega\in \Omega^1$ and $\eps'=\pm1$. For an `even' spectral triple, we also require a bimodule map $\gamma:\CS\to \CS$ such that 
\begin{equation}\label{gamma} \gamma^2=\id,\quad \gamma\circ\la= - \la (\id\tens\gamma),\quad \CJ\circ\gamma= \eps'' \gamma\circ \CJ\end{equation}
for a third sign $\eps''=\pm1$. The pattern of signs then determines the `spectral triple' dimension mod 8 in Connes' theory\cite{Con,ConMar}. We also ask that  $\gamma:\CS\to \CS$ intertwines $\nabla_{\CS}$ on each side,
\begin{equation}\label{gammanabla} \nabla_{\CS}\circ\gamma= (\id\tens\gamma)\nabla_{\CS}.\end{equation}

\subsection{}\label{sec25} The geometric approach \cite{BegMa:spe,BegMa} proposes two further conditions not required for a spectral triple but natural for the geometry. The first is to be compatible with 
a QLC on $\Omega^1$ in such a way that $\la$ intertwines the tensor product connection and the one on $\CS$,
\begin{equation}\label{qlcnablaS} \nabla_{\CS}\circ\la= (\id\tens\la)\nabla_{\Omega^1\tens_A\CS},\end{equation}
i.e. covariance $\doublenabla(\la)=0$. The second is compatibility with $\Omega^2$ in the sense that $\la$ extends to a left module map $\Omega^2\tens_A\CS\to \CS$, for example according to 
\begin{equation}\label{Omega2la} \varphi(\omega\la(\eta\la s))=\kappa (\omega,\eta)s+(\omega\wedge\eta)\la s\end{equation}
 for all $\omega,\eta\in \Omega^1$ and $s\in \CS$, for some constant $\kappa$ and some invertible left module map $\varphi:\CS\to \CS$. We will only need the nicest case where $\varphi=\id$ and $\kappa=1$. 

\subsection{}\label{sec26}  It is shown in \cite{BegMa:spe} that (\ref{J})--(\ref{gammanabla})  imply the algebraic part of the axioms of a Connes spectral triple\cite{Con}, i.e.  {\em other} than those that relate to the Hilbert space inner product. For the latter, we further require a hermitian inner product on $\CS$ so as to be able to complete it to a Hilbert space on which $A$ is represented as a $*$-algebra. We then require that $\imath\dirac$ (in our conventions) and $\gamma$ (if it exists) are hermitian and $\CJ$ is an antilinear isometry in the sense $\<\CJ\psi,\CJ\phi\>=\<\phi,\psi\>$ for all $\phi,\psi\in \CH$. Details on this division of Connes' axioms and the precise construction from our geometric data are in \cite[Sec.~8.5]{BegMa}. 

\subsection{}\label{secsetup} We will be interested in the case where $A$ has trivial centre, $\Omega^1$ has a central basis $\{s^i\}$, say, and $\CS$ has a central basis $\{e^\alpha\}$, say. We later apply this with $i=1,2,3$ and  $\alpha=1,2$ to the fuzzy sphere.  Some general analysis for quantum Riemannian geometry on $\Omega^1$ in this parallelizable case is in \cite[Ex.~8.2 ]{BegMa} and we extend this now to include the above constructions on $\CS$. 

In this case,  we can write the  `Clifford action' $\la:\Omega^1\tens_A \CS\to \CS$ as given by a matrix $C^{i\alpha}{}_\beta\in A$ according to
\begin{equation}\label{Cia} s^i\la e^\alpha= C^{i\alpha}{}_\beta e^\beta .\end{equation}
Centrality of the basis and the trivial centre forces the coefficients to be  $C^{i\alpha}{}_\beta\in\C$. We will also consider them as a collection of matrices $(C^i){}^\alpha{}_\beta=C^{i\alpha}{}_\beta$.  Next, we define 
\[ \CJ(a e^\alpha)= a^*J^\alpha{}_\beta e^\beta\]
 and centrality of the basis and the trivial centre mean that (\ref{J}) requires the $J^\alpha{}_\beta\in \C$ and then  
\begin{equation}\label{JJ}J^\alpha{}_\beta{}^* J^\beta{}_\gamma=\eps\delta^\alpha{}_\gamma;\quad \overline{J}J=\eps\id,\quad \eps=\pm1,\end{equation}
where we also write the equation in a compact form for a matix $J$. The overline denotes complex conjugation of the entries.

We similarly write $\nabla_{\CS}=S^\alpha{}_{i\beta}s^i\tens e^\beta$ for a collection of matrices $(S_i){}^\alpha{}_\beta=S^\alpha{}_{i}{}_\beta$,  which we suppose have entries in $\C$.  This is equivalent to supposing that $\sigma_{\CS}$ is the flip on the bases. In this case,  the reality condition is
\begin{equation}\label{SJ}
            S^\alpha{}_{i \beta}{}^* J^\beta{}_{\gamma}=J^\alpha{}_\beta S^\beta{}_{i \gamma} ;\quad \overline{S_i}J=JS_i,
         \end{equation}
while the second part of (\ref{JnablaS}), for compatibility of $\CJ$ with the Clifford action, translates to 
\begin{equation}\label{CJ}
            C^{i\alpha}{}_{\beta}{}^* J^\beta{}_{\gamma} = \epsilon' J^\alpha{}_{\beta}C^{i \beta}{}_{\gamma} ;\quad \overline{C^i}J=\eps'JC^i,\quad \eps'=\pm 1.
          \end{equation} 

Finally, for an  `even' spectral triple, we also set $\gamma (e^\alpha) = \gamma^\alpha{}_{\beta} e^\beta$, here again with $\gamma^\alpha{}_{\beta} \in \C$ under our central basis and centre assumptions. Then  the three parts of (\ref{gamma}) are respectively
\begin{align}\label{gamsq}
          \gamma^\alpha{}_{\beta} \gamma^\beta{}_{\nu} &= \delta^{\alpha}{}_{\nu};\quad\quad\qquad \quad\gamma^2=\id,
         \\ \label{Cgam} 
          C^{i\alpha}{}_{\beta}\gamma^{\beta}{}_{\nu}&=-\gamma^\alpha{}_{\beta} C^{i\beta}{}_{\nu};\quad\qquad \{C^i,\gamma\}=0,          \\ \label{Jgam}
          \gamma^\alpha{}_{\beta}{}^* J^\beta{}_{\nu} &= \epsilon'' J^\alpha{}_{\beta} \gamma^\beta{}_{\nu};\quad\qquad \overline{\gamma}J=\eps''J \gamma,\quad\eps''=\pm1.
          \end{align}
Note that the full set of conditions for $J$ is invariant under multiplication by a phase. We also need compatibility (\ref{gammanabla}) with the connection,  
\begin{equation}\label{Sgam}
          S^\alpha{}_{i\beta}\gamma^{\beta}{}_{\nu}=\gamma^\alpha{}_{\beta} S^\beta{}_{i\nu};\quad [S_i,\gamma]=0.
         \end{equation}

 Beyond the local data for  a spectral triple, we have the further geometric requirements in Section~\ref{sec25}. Here the covariance  $\doublenabla(\la)=0$  of the Clifford action expressed in (\ref{qlcnablaS}) amounts to 
 \begin{equation}\label{CSGam}
          C^{i \alpha}{}_{\beta} S^\beta{}_{j\nu}-S^\alpha{}_{j\beta}C^{i \beta}{}_{\nu} = -\tfrac{1}{2} \tchrs{i}{jk} C^{k\alpha}{}_{\nu};\quad [C^i,S_j]=-\tfrac{1}{2} \tchrs{i}{jk}C^k,
          \end{equation}
 where we suppose in keeping with the metric that $\nabla s^i=- {1\over 2}\tchrs{i}{jk}s^j\tens s^k$ has constant coefficients.  Finally, we suppose that $s^i$ form a Grassmann algebra as a basis over $A$ of $\Omega^2$, in which case compatibility (\ref{Omega2la}) with the Clifford action  becomes
 \begin{equation}\label{CC}
       C^{i\alpha}{}_{\nu} C^{j\nu}{}_{\beta}+ C^{j\alpha}{}_{\nu} C^{i\nu}{}_{\beta} = 2g^{ij}\delta^\alpha{}_\gamma;\quad \{C^i,C^j\} =2g^{ij}\id        \end{equation}       
for the nicest case where $\varphi=\id$ and $\kappa=1$. This just says that $C^i$ represent a usual Clifford algebra. 

\section{Construction of the Dirac operator}\label{secround}

Here we consider the construction of the Dirac operator on $A=\C_\lambda[S^2]$ in the geometric form above, initially for the case of the round metric $g_{ij}=\delta_{ij}$ and quantum Levi-Civita connection $\Gamma^i{}_{jk}=\eps_{ijk}$. We show that this leads to a natural rotationally invariant Dirac operator meeting Connes' axioms of a spectral triple on completion of the spinor bundle $\CS$ to a Hilbert space. Clearly, the standard Clifford structure of $\R^3$ (so $C^i=\sigma^i$, the Pauli matrices) is then the obvious solution of (\ref{CC}) and we start with this case. 

\begin{lemma}\label{lemmoduli} For the round metric quantum geometry and $\la$ given by the standard Clifford structure of $\R^3$, there is a 4-parameter moduli of $\nabla_S,\CJ$ obeying the remaining conditions (\ref{JJ})--(\ref{CSGam}), with $\eps=\eps'=-1$. Here,
\[\CJ(a e^1)=a^* qe^2,\quad \CJ(ae^2)=-a^*qe^1,\quad s^i\la e^\alpha=\sigma^i{}^\alpha{}_\beta e^\beta,\]
\[\nabla_\CS e^\alpha={\imath\over 4}\sigma^i{}^\alpha{}_\beta s^i\tens e^\beta+d_i s^i\tens e^\alpha,\quad \dirac (\psi_\alpha e^\alpha)=( (d_i+\del_i) \psi_\beta\sigma^i{}^\beta{}_\alpha+{3 \imath\over 4} \psi_\alpha ) e^\alpha \]
with parameters $q$ a phase and $d_i\in\R$. We give $\dirac$ on a spinor  $\psi_\alpha e^\alpha\in\CS$ with coefficients $\psi_\alpha\in A$.    
\end{lemma}
\proof The choice of $C^i$ already solves the `Clifford action' condition (\ref{CC}). We set
\[ S^\alpha{}_{i\beta}={\imath\over 4} C^{i\alpha}{}_\beta+d_i \delta^\alpha{}_\beta,\quad d_i\in \C;\quad S_i={\imath\over4}C^i+d_i\id,\]
which automatically solves (\ref{CSGam}).  One can check that this is the most general solution. In our analysis, this gives a bimodule connection with $\sigma_\CS$ the flip map. Next, (\ref{CJ}) forces $\eps'=-1$ and
\[ J=q\begin{pmatrix} 0 & 1\\ -1& 0\end{pmatrix},\quad q\in \C,\]
where only $\sigma_2$ and hence the matrix $C^2{}^\alpha{}_\beta$ has complex (imaginary) entries. Now imposing (\ref{JJ})  requires
\[ \eps=-1,\quad |q|=1.\]
Finally, (\ref{SJ})  holds automatically provided $d_i$ are real since any multiple $d_i\id$ commutes with $J$, and $C^i$ already obey (\ref{CJ}) with $\eps'=-1$. 
One can show that there is no nonzero matrix $\gamma$ for which (\ref{Cgam}) holds so we drop the group (\ref{gamsq})-(\ref{Sgam}).    Next, we feed our matrix solutions into the general construction to get the Clifford action, $\CJ$,  $\nabla_\CS$ and hence $\dirac$ on the $e^\alpha$ basis as stated.  Note that the signs $\eps,\eps'$ look like they fit in  Connes' axioms for KO dimension $n=5$ mod 8 but we will see that the actual spectral triple has $\imath\dirac$ not $\dirac$ which changes the analysis. 

In terms of the coefficients of the spinors, one can write
\begin{equation}\label{diraccomponents} (\dirac\psi)_\alpha=(d_i+\del_i)\psi_\beta\sigma^{i\beta}{}_\alpha+{3\imath\over 4}\psi_\alpha,\end{equation}
where $\psi_\alpha$ are (in our conventions) a co-spinor (a row vector in spinor space). If we use column vector notation then this would appear as $\dirac= \sigma^i{}^T(d_i+\del_i)+{3\imath \over 4}$. \endproof

The theory in \cite{BegMa:spe} implies that so far we obey Connes' axioms for a spectral triple other than those involving the Hilbert space inner product. To address the latter, we need a notion of integral $\int:A\to \C$, for which, we propose to use the noncommutative spherical harmonic expansion as in \cite{ArgMa},
\[ A=\oplus_{l\in \N\cup\{0\}} A_l;\quad A_l=\{f_{i_1\cdots i_l}x^{i_1}\cdots x^{i_l}\ |\ f\ {\rm totally\ symmetric\ and\ trace\  free}\}\]
A similar expansion was used for the reduced  $M_n(\C)$ fuzzy spheres in \cite{Mad} with $l=0,\cdots,n-1$. We define $\int a$ to be the $A_0$ component of $a$ in this expansion. 

\begin{lemma}\label{lemround} $\int$ on $A=\C_\lambda[S^2]$ is rotationally invariant and obeys
\[ \int \del_i a=0,\quad \int ab= \int ba,\quad \int a^*=\overline{\int a}\]
for all $a,b\in A$. We also have $\int a^*a\ge 0$ with equality if and only if $a=0$. 
\end{lemma}
\proof  By definition, $\int$ is clearly rotationally invariant. As $\del_i={1\over 2\imath\lambda_p}[x_i,\ ]$ acts as orbital angular momentum, rotational invariance implies $\int\del_i a=0$. It also implies that $\int x_i a=\int a x_i$, and iterating this, that $\int ab=\int ba$ for all $a,b$. The behaviour under complex conjugation is also clear given that $x_i^*=x_i$. Positivity is less clear, but one way is to see it is to suppose a completion of $\C_\lambda[S^2]$ to a $C^*$-algebra  and define
\[ \int a=\int_{SU_2} \extd g\, {\rm Rot}_g(a)\]
where the Haar measure on $SU_2$ is normalised so that $\int 1=1$ and  ${\rm Rot}$ is the rotation action on $\C_\lambda[S^2]$ by $SU_2$. This picks out the $l=0$ part. Its image lies in the centre of $\C_\lambda[S^2]$ which we identify with $\C.1$.  From this point of view,  rotational invariance and the behaviour under complex conjugation is again clear, but so is positivity  as $\int\extd g\, {\rm Rot}_g(a^*a)=\int\extd g\,  {\rm Rot}_g(a)^*{\rm Rot}_g(a)$ is a convex linear combination of  positive operators. We illustrate these properties in low degree. Thus,
\begin{align*}
\del_i(x_j x_k)&=\eps_{ijm}x_m x_k+\eps_{ikm}x_jx_m= (\eps_{ijl}\delta_{km}+\eps_{ikm}\delta_{jl})x^lx^m\\
&={1\over 2}(\eps_{ijl}\delta_{km}+\eps_{ikm}\delta_{jl}+\eps_{ijm}\delta_{kl}+\eps_{ikl}\delta_{jm})x^lx^m+{1\over 2}(\eps_{ijl}\delta_{km}+\eps_{ikm}\delta_{jl})[x^l,x^m]\\
&={1\over 2}(\eps_{ijl}\delta_{km}+\eps_{ikm}\delta_{jl}+\eps_{ijm}\delta_{kl}+\eps_{ikl}\delta_{jm})x^lx^m+{\imath\lambda_p}(\eps_{ijl}\delta_{km}+\eps_{ikm}\delta_{jl})\eps_{lmp}x^p\\
&={1\over 2}(\eps_{ijl}\delta_{km}+\eps_{ikm}\delta_{jl}+\eps_{ijm}\delta_{kl}+\eps_{ikl}\delta_{jm})x^lx^m+{\imath\lambda_p}(\delta_{ik}\delta_{jp}-\delta_{ij}\delta_{kp})x^p
\end{align*}
where the first tensor is symmetric traceless in the $l,m$ indices. Hence there is no $l=0$ component and the integral is zero. By contrast, 
\begin{equation}\label{intxij} x_ix_j={1\over 3}(1-\lambda_p^2)\delta_{ij}+\imath\lambda_p\eps_{ijk}x^k+ {1\over 2}(\delta_{ik}\delta_{jl}+\delta_{jk}\delta_{il}-{2\over 3}\delta_{ij}\delta_{kl} )x^kx^l\end{equation}
is an expansion into $l=0,1,2$ components using the relations of the algebra, so $\int x_ix_j={1\over 3}(1-\lambda_p^2)\delta_{ij}$, which is positive for $i=j$ and implies positivity on $A_1$. 
\endproof

We now define on $\CS$ the sesquilinear inner product
\begin{equation}\label{Sinnerprod}\<\psi_\alpha e^\alpha, \phi_\beta e^\beta\>=\int \psi_\alpha{}^*\phi_\alpha,\end{equation}
which obeys
\[ \overline{\<\psi_\alpha e^\alpha, \phi_\beta e^\beta\>}=\int \phi_\alpha{}^*\psi_\alpha= \<\phi_\beta e^\beta, \psi_\alpha e^\alpha\>\]
as required, as well as positivity. This therefore completes to a Hilbert space 
\[ \CH=\overline{\CS}=\{\psi=\sum_{l=0}^{\infty} \psi_l,\quad \psi_l\in\CS_l\ |\ \sum_{l=0}^\infty ||\psi_l|^2 <\infty\},\]
where each $\CS_l=A_l\oplus A_l$ is a finite-dimensional sub-Hilbert space and these subspaces are mutually orthogonal since the product of elements of different orbital angular momentum has no component in $l=0$. We define $\CJ,\dirac$ by their values on $\CS$.  

\begin{proposition}\label{dirround} In Lemma~\ref{lemround}, $\CJ$ is an antilinear isometry and if $d_i=0$ then $\imath\dirac$ is hermitian. In this case, $(A,\CH,\CJ,\imath\dirac)$ is a geometrically realised KO-dimension $n=3$ Connes spectral triple on the fuzzy sphere.
\end{proposition}
\proof  We first check $\CJ$, 
\begin{align*}\<\CJ(\phi_\beta e^\beta),\CJ(\psi_\alpha e^\alpha)\>&=\<\phi_\beta{}^*J^\beta{}_\gamma e^\gamma, \psi_\alpha{}^*J^\alpha{}_\delta e^\delta\>=\int \overline{J^\beta{}_\gamma}\phi_\beta \psi_\alpha{}^* J^\alpha{}_\gamma\\
&=-\int \overline{J^\beta{}_\gamma}J^\gamma{}_\alpha \phi_\beta \psi_\alpha{}^*=\int \phi_\alpha \psi_\alpha{}^*=\int  \psi_\alpha{}^* \phi_\alpha=\<\psi_\alpha e^\alpha, \phi_\beta e^\beta\>\end{align*}
as required, using that $J$ in our case is antisymmetric for the 3rd equality. For $\dirac$,
\begin{align*} \<\dirac(\psi_\alpha e^\alpha),\phi_\beta e^\beta\>&=\int (-{3\imath\over 4}\psi_\alpha{}^*+ (d_i +\del_i)\psi_\beta{}^*\overline{\sigma^{i\beta}{}_\alpha})\phi_\alpha\\
&=-\int \psi_\alpha{}^*({3\imath\over 4}\phi_\alpha+ (-d_i+\del_i)\phi_\beta\sigma^{i\beta}{}_\alpha=-\<\psi_\alpha e^\alpha,\dirac(\phi_\beta e^\beta)\>\end{align*}
{\em provided} $d_i=0$. We used for the 2nd equality that $\sigma^i$ are hermitian and integration by parts $\int (\del_i\psi_\alpha{}^*)\phi_\alpha=\int\del_i (\psi_\alpha{}^*\phi_\alpha)-\int \psi_\alpha{}^* \del_i\phi_\alpha=-\int \psi_\alpha{}^* \del_i\phi_\alpha$ given properties of the integral and that $\del_i$ is a derivation because the $s^i$ are central. We also used that $\del_i$ commutes with $*$ because $s^i{}^*=s^i$.  Thus, $\imath\dirac$ is symmetric with respect to the inner product. We will see shortly below that it is diagonalisable and preserves $\CS_l$. It is therefore bounded and self-adjoint on each subspace, from which it is clear that the domain of $\dirac$ and its adjoint are the same, so that $\imath\dirac$ is self-adjoint. Finally, because $\CJ$ is antilinear, it means that $\imath\dirac$ commutes with $\CJ$, not anticommutes as found for $\dirac$. Allowing for this, the sign relevant for $\imath\dirac$ is  $\eps'=1$, which means KO dimension 3 mod 8. \endproof

Looking at the geometric derivation, it is striking that we more naturally obtain $\dirac$ antihermitian not hermtian (albeit easily fixed with an $\imath$). The normalisations of $C^i,S_j$ are determined by the geometric conditions (\ref{CSGam})-(\ref{CC}) which we optionally imposed but which are natural and indeed fit with the geometric Dirac operator being a derivative and not $\imath$ times a derivative. On the other hand, we see that being antihermitian eliminates the unexpected $d_i$ part of the form of $\dirac$. The $\imath\over 4$ reflects what classically would be the constant curvature in the Lichnerowicz formula. The quantum Lichnerowicz  formula is in \cite[Prop.~8.45]{BegMa} and relates  the Laplace-Beltrami operator on $\CS$  defined by $\Delta_S=((\ ,\ )\tens\id)\nabla_{\Omega^1\tens S}\nabla_S$ and the curvature of $\nabla_\CS$ defined by $R_{\CS}=(\extd\tens\id-\id\wedge\nabla_{\CS})\nabla_{\CS}$.
 
 \begin{proposition} For $\dirac$ in Proposition~\ref{dirround}, we have $\dirac^2=\Delta_\CS+\la\circ R_\CS$,  
 where
 \[ \Delta_S(\psi_\alpha e^\alpha)=(\del_i\del_i \psi_\alpha+ {\imath\over 2}\del_i\psi_\beta \sigma^{i\beta}{}_{\alpha} -{3\over 16} \psi_\alpha) e^\alpha,\quad  \la\circ R_\CS=-{3\over 8}\id\]
 \end{proposition}
 \proof This form of $\dirac^2$ is from \cite[Prop.~8.45]{BegMa} in the general theory for geometrically realised Dirac operators, with $\varphi=\id$ and $\kappa=1$ there. We use that our $\nabla$ on $\Omega^1$ has zero torsion and that $\doublenabla(\la)=0$ as the optional  intertwining  condition that we required. We now compute the two parts explicitly, starting with  
 \begin{align*} \Delta_\CS(\psi_\alpha e^\alpha)&=((\ ,\ )\tens\id)\nabla_{\Omega^1\tens\CS}((\del_i \psi_\alpha+ \psi_\beta S^\beta{}_{i\alpha})s^i\tens e^\alpha)\\
 &=  \del_j(\del_i \psi_\alpha+ \psi_\beta S^\beta{}_{i\alpha})(s^j, s^i) e^\alpha+(\del_i \psi_\alpha+ \psi_\beta S^\beta{}_{i\alpha})((\ ,\ )\tens\id)\nabla_{\Omega^1\tens\CS}(s^i\tens e^\alpha)\\
 &=\del_i(\del_i \psi_\alpha+ \psi_\beta S^\beta{}_{i\alpha}) e^\alpha+(\del_i \psi_\alpha+ \psi_\beta S^\beta{}_{i\alpha})(-{1\over 2}\eps_{ijk}(s^j,s^k)e^\alpha+(s^j,s^i)S^\alpha{}_{j\gamma}e^\gamma)\\
 &=(\del_i\del_i \psi_\alpha+ {\imath\over 2}\del_i\psi_\beta \sigma^{i\beta}{}_{\alpha})e^\alpha+ \psi_\beta({\imath\over 4})^2\sigma^{i\beta}{}_\gamma\sigma^{i\gamma}{}_\alpha e^\alpha, 
 \end{align*}
 which simplifies as stated on using the Pauli matrix identity $\sigma^i\sigma^j=\delta_{ij}+\imath\eps_{ijk}\sigma^k$. Meanwhile, $\la R_\CS$ is a left module map so it suffices to give it on the basis,
 \begin{align*} \la R_\CS(e^\alpha)&= \la (\extd \tens\id -\id\wedge\nabla_\CS)(s^i\tens S^{\alpha}{}_{i\beta}e^\beta)\\ &=S^{\alpha}{}_{i\beta}(-{1\over 2}\eps_{ijk}s^j\wedge s^k\la e^\beta- s^i\wedge s^j \la S^\beta{}_{j\gamma}e^\gamma)\\
 &=-{1\over 2} \eps_{ijk}S^{\alpha}{}_{i\beta}C^{k\beta}{}_\gamma C^{j\gamma}{}_\delta e^\delta-S^{\alpha}{}_{i\beta}S^\beta{}_{j\gamma}C^{j\gamma}{}_\delta C^{i\delta}{}_\eta e^\eta+S^{\alpha}{}_{i\beta}S^\beta{}_{j\gamma}(s^i,s^j)e^\gamma \\
 &=(-{\imath\over 8}\eps_{ijk}\sigma^i\sigma^k\sigma^j+{1\over 16}\sigma^i\sigma^j\sigma^j\sigma^i-{1\over 16}\sigma^i\sigma^i)^\alpha{}_\beta e^\beta=(- {3\over 4}+{9\over 16}-{3\over 16})e^\alpha=-{3\over 8}e^\alpha,\end{align*}
where we used $s^i\wedge s^j\la e^\gamma=s^i\la(s^j\la e^\gamma)-(s^i,s^j)e^\gamma$ for the third equality. It is not necessary, but as a check one can also verify the result by computing $\dirac^2$ directly. Thus, 
 \begin{align*}\dirac^2(\psi_\alpha e^\alpha)&=\dirac\left(({3\imath\over 4}\psi_\alpha+ \del_i \psi_\beta \sigma^{i\beta}{}_\alpha)e^\alpha\right)\\
 &=\left({3\imath\over 4}({3\imath\over 4}\psi_\alpha+\del_i\psi_\beta \sigma^{i\beta}{}_\alpha)+\del_j({3\imath\over 4}\psi_\beta+\del_i\psi_\gamma\sigma^{i\gamma}{}_\beta)\sigma^{j\beta}{}_\alpha\right)e^\alpha\\
 &=\left(-{9\over 16}\psi_\alpha+ {3\imath\over 2}\del_i \psi_\beta\sigma^{i\beta}{}_\alpha+(\del_j\del_i\psi_\gamma)(\sigma^i\sigma^j)^\gamma{}_\alpha\right)e^\alpha\\
 &=\left(-{9\over 16}\psi_\alpha+{\imath\over 2}\del_i\psi_\beta\sigma^{i\beta}{}_\alpha+\del_i\del_i \psi_\alpha\right)e^\alpha
 \end{align*}
 using the Pauli matrix identity and
 \[ \eps_{ijk}\del_i\del_j=\del_k\]
 as follows from $\del_i={1\over 2\imath\lambda_p}[x_i, ]$, the Jacobi identity and the $[x_i,x_j]$ commutation relations. This therefore agrees with the result from \cite[Prop.~8.45]{BegMa}. \endproof 
 
\begin{proposition} $\imath\dirac$ in Proposition~\ref{dirround} is diagonalisable with eigenvalues $\lambda_{l,\pm}$ of multiplicity $n_{l,\pm}$,  
\[  \lambda_{l,\pm}=-{1\over 4}\pm (l+{1\over 2});\quad n_{l,\pm}=2l+1\mp 1\] 
according to a $\pm$ sign  and orbital angular momentum $l=0,1,2,\cdots$. Here $l=0$ only has $\lambda_{0,-}$. \end{proposition}
\proof Here $\del.\del$ commutes with $\dirac$ so they can be simultaneously diagonalised. Hence it suffices to diagonalise $\dirac$ on the subspaces $\CS_l$ of fixed orbital angular momentum $l$, where $\del_i\del_i=-l(l+1)$ by \cite[Prop.~4.3]{ArgMa} for each spinor component. From $\dirac^2$ above, we see that
\[ -(\imath \dirac)^2-{\imath \dirac\over 2}=\sum_i\del_i^2 -{3\over 16},\]
which gives the possible eigenvalues of $\imath\dirac$ for fixed $l$. For $l=0$, we obviously only have $-3/4$ but for higher spin one can check that both eigenvalues occur. For example, for $l=1$ there are (by direct calculation) 2 eigenvectors of eigenvalue $5/4$ and 4 eigenvectors of eigenvalue $-7/4$, namely (in a 2-vector notation for spinor space)
\[ \lambda_{1,+}:\quad \begin{pmatrix}x_1+\imath x_2\\ -x_3\end{pmatrix},\ \begin{pmatrix}x_3\\ x_1-\imath x_2\end{pmatrix};\quad \lambda_{1,-}:\quad \begin{pmatrix}x_1\\ x_3\end{pmatrix},\ \begin{pmatrix}x_3\\ -x_1\end{pmatrix},\ \begin{pmatrix}0\\ x_1+\imath x_2\end{pmatrix},\  \begin{pmatrix}x_1-\imath x_2\\0\end{pmatrix}\] 
forming multiplets of total spin 1/2 and 3/2 respectively. In general, $\CS_l= A_l\tens \C^2$ is a direct sum of a total spin $l-{1\over 2}$ irreducible of $SU_2$ and a total spin $l+{1\over 2}$ irreducible. These add up to the dimension $2(2l+1)$ of $\CS_l$ since each component $A_l\subset \C_\lambda[S^2]$ is a $2l+1$-dimensional irreducible represention. Hence, this totally diagonalises $\dirac$. This also completes the proof of Proposition~\ref{dirround}.  \endproof

This spectrum differs from the usual spectrum $\lambda_{l,\pm}=\pm(l+1)$ with multiplicity $2(l+1)$ for a classical sphere $S^2$, but has in common the key feature that the spectrum is discrete with no zero modes. Indeed, adding $-1/4$ to our spectrum gives the same set $\pm 1,\pm 2, \cdots$  etc. with the same multiplicities as classically, so the difference can be viewed as this constant shift, although the underlying oribital angular momentum is a little different in terms of how it contributes to the two series.  This amounts to a relatively contained impact of the fuzzy sphere's extra non-classical cotangent direction. The latter was seen in \cite{LirMa} to also change the Ricci curvature of $-{3/4}$ in our conventions (which would be 3/2 in usual conventions, compared to the classical value of 2 for a unit sphere) and to have other more drastic consequences in \cite{ArgMa}.

So far, we have adopted the standard Clifford structure given by the Pauli matrices  and proceeded from there. For completeness, we show that any other solutions of our spectral triple equations are equivalent to the one above. 

\begin{proposition}\label{dirU} Equations (\ref{JJ})--(\ref{CC}) for a geometric spectral triple and $g_{ij}=\delta_{ij}$ have the most general solution
\[ C^i=U\sigma_iU^{-1},\quad S_i={\imath\over 4} C^i+d_i\id,\quad J=q\overline{U}\begin{pmatrix}0&1\\-1&0\end{pmatrix}U^{-1};\quad |q|=1,\quad d_i\in\R,\quad \eps=\eps'=-1,\]
with no solution for $\gamma$ (so an odd spectral triple), $q$ a phase parameter and $U\in SL_2(\C)$. Moreover, $\imath \dirac$ is Hermitian if and only if $d_i=0$ and $U$ is unitary, in which case $\CJ$ is an antilinear isometry. Hence, any spectral triple constructed by our method is unitarily equivalent to the one in Proposition~\ref{dirround}, where $U=\id$.  \end{proposition}
\proof It is known from the theory of Clifford algebras that any solution of (\ref{CC})  and $g_{ij}=\delta_{ij}$ with be conjugated to our standard Pauli-matrix one, but for completeness one can also show this, as follows. Thus, the $C^i$ square to the identity and the Clifford algebra relations are invariant under conjugation of the $C^i$. But a Jordan 2-block can never square to the identity. Hence, without loss of generality, we can suppose that $C^3$ is either $\pm$ the identity or $\sigma^3$. Now solving for the remaining Clifford algebra relations, one finds that only $C^3=\sigma^3$ works and that there there is a 1-parameter family of solutions for $C^1,C^2$. These are, however, still conjugate to our previous solution given by the $\sigma^i$, by a diagonal $U$ that does not change $C^3$.

Using this observation, we can choose a basis rendering any $C^i$ to be our standard one. As (\ref{CSGam}) is linear and invariant under conjugation, we can solve for $S_j$ in the preferred basis. It follows that $S_i={\imath \over 4}C^i+d_i\id$ in any basis for some $d_i\in \C$ and that it has the form stated. It similarly follows that there are no solutions for $\gamma,\eps''$ by looking at  (\ref{Cgam}). In this case, comparing (\ref{SJ}) and  (\ref{CJ}) forces $\eps'=-1$ and $d_i$ real. 

It remains to solve for $\eps',J$ for fixed invertible $U$ of (say) determinant 1 (as multiplying by a scale does not change $C^i,S^i$). Following the same format as in the proof of Lemma~\ref{lemmoduli} but now for  $C^i,S^i$ defined by $U$, we look for solutions of (\ref{CJ}) and find that this needs $\eps'=-1$ and results in a particular 1-parameter family of $J$. Then requiring (\ref{JJ}) forces $\eps=-1$ as before but now  
\[ J=
q\begin{pmatrix}
 - u_{21} \overline{u^1{}_1}-u^2{}_2 \overline{u^1{}_2} &  | u^1{}_1| ^2+| u^1{}_2| ^2 \\
 - |u^2{}_1|^2-|u^2{}_2 |^2 &  u^1{}_1 \overline{u^2{}_1}+u^1{}_2 \overline{u^2{}_2}
\end{pmatrix};\quad |q|=1\]
for a phase parameter $q$ and $U=(u^\alpha{}_\beta)$, which we write as stated. Clearly, $U=\id$ recovers our previous  spectral triple data. By our analysis, this is also the most general solution of (\ref{JJ})--(\ref{CC})  under our assumptions. 

Next, for a spectral triple, we ask for $\imath\dirac$ to be symmetric. Looking  in the proof of Proposition~\ref{dirround}, we need that $d_i=0$ and that $C^i{}^\dagger=C^i$ or $U\sigma^iU^{-1}=(U\sigma^i U^{-1})^\dagger$ for all $i$, which amounts to $[\sigma^i,U^\dagger U]=0$ and hence that $U^\dagger U$ is a multiple of the identity. From the determinant, this multiple is $\pm1$. But ${\rm Tr}(U^\dagger U)=\sum_{\alpha,\beta}|u^\alpha{}_\beta|^2>0$ excludes the $-1$ case. Hence we are forced to $U$ unitary. In this case, one can show that $J$ remains antisymmetric, e.g.  in the expression stated for $J$, one has
\[ |u^1{}_1|^2+|u^1{}_2|^2-|u^2{}_1|^2-|u^2{}_2|^2=2(|u^1{}_2|^2-|u^2{}_1|^2)=0\]
for a unitary matrix. This antisymmetry and (\ref{JJ}) were what we needed in the  proof of Proposition~\ref{dirround}, so this goes through as before and $\CJ$ is an antilinear isometry. 

Finally, with $d_i=0$, since $\dirac$ only depends on $C^i$, if $\psi_\alpha e^\alpha$ is an eigenfunction of our original Dirac operator, then clearly $\psi'_\alpha=\psi_\beta U^{-1}{}^\beta{}_\alpha$ is an eigenfunction for our new Dirac operator with conjugated $C^i$. This has the same inner product (\ref{Sinnerprod}) when $U$ is unitary, so the eigenvalues are  unchanged and the self-adjointness of $\imath\dirac$ also holds. Hence, this spectral triple is unitarily equivalent to our original one given by $U=\id$.  \endproof

It remains to consider what happens in our approach when $g_{ij}$ is not proportional to the round metric. Note that a general quantum metric has the form 
\begin{equation}\label{gO} g_{ij}=O_{ik}O_{jl}\delta_{kl}\lambda_k,\quad g^{ij}=O_{ik}O_{jl}\delta_{kl}\lambda_k^{-1}\end{equation}
for some  $O\in SO(3)$ and some nonzero real eigenvalues $\lambda_k$. In this case the Clifford algebra relation (\ref{CC})  is  solved by
\begin{equation}\label{CO}  C^i=O_{ij}\sigma^j\lambda_j^{-{1\over 2}}.\end{equation}
This is not the only solution but, as before, we can conjugate any solution by $U$ in spinor space to put it in this canonical form. We focus on the Euclidean signature where $\lambda_i>0$.

\begin{proposition}\label{propdirg} Following the same construction (\ref{JJ})--(\ref{CC}) as above but with general $g_{ij}$ of Euclidean signature, leads to a Dirac operator
\[  \dirac(\psi_\alpha e^\alpha)=(\del_i\psi_\alpha C^{i\alpha}{}_\beta + d_i \psi_\alpha\sigma^{i\alpha}{}_\beta)e^\beta+{\imath\over 4}{{\rm Tr}(g)\over\sqrt{\det(g)}}\psi_\alpha e^\alpha \]
for some real constants $d_i$. This forms a spectral triple with $\imath\dirac$ self adjoint if and only if $d_i=0$. 
\end{proposition}
 \proof We have already fixed our solution to (\ref{CC}) as (\ref{CO}). We next compute for the quantum Levi-Civita connection in (\ref{QLC})  that
 \begin{align*}\Gamma^i{}_{jk}O_{kl}&=O_{im}{\delta_{mq}\over\lambda_m}O_{pq}(2\eps_{pkr}O_{rs}O_{jn}\delta_{sn}\lambda_n+{\rm Tr}(g)\eps_{pjk} )O_{kl}\\
&=O_{im}O_{jn}{\delta_{mq}\over\lambda_m}\left(2\eps_{qls}\delta_{sn}\lambda_n+ {\rm Tr(g)}\eps_{qnl}\right)\\
&=O_{im}O_{jn}{\eps_{mnl}\over\lambda_m}\left(-2\lambda_n+{\rm Tr}(g)\right)
 \end{align*}
 using $\det(O)=1$ so that $\eps_{ijk}$ is invariant. Using this expression in (\ref{CSGam}) together with (\ref{CO}) and  $S_i=O_{ij}\varsigma_j$ for some spinor matrices $\varsigma_i$,  this becomes for all $i,j$ (and summing over $k$), 
 \[ [\sigma^i,\varsigma_j]=\eps_{ijk}\sigma^k c_{ijk};\quad c_{ijk}=-{1\over 2\sqrt{\lambda_i\lambda_k}} (-2\lambda_j+\sum_m\lambda_m).\]

Now suppose that $\varsigma_i=\beta_{ij}\sigma^j$ for some coefficients $\beta_{ij}$ (plus a possible multiple of the identity which does not change the analysis). Then this requires for all $i,j,k$ (summing over $m$), 
\[ 2\imath\beta_{jm}\eps_{imk}= \eps_{ijk} c_{ijk}.\]
Setting $k=j\ne i$, we learn that $\beta$ is diagonal, while from the six distinct values of $ijk$ we learn that $\beta_{jj}=c_{ijk}/(2\imath)$ independently of the order of the other two indices. Putting all this together, (\ref{CSGam}) is solved by
\begin{equation}\label{Sigen} S_i=O_{ij}({\imath\sigma^j\over 4}\mu_j+d_j){\sqrt{\lambda_j}};\quad \mu_j={1\over\sqrt{ \lambda_1\lambda_2\lambda_3}}(-2\lambda_j+\sum_m\lambda_m)\end{equation}
for some constants $d_i$. The equations for $J$ are the same as when $O=\id$ since they apply for each $i$ and $O_{ij},\mu_i$ are real. Hence we solve (\ref{JJ})--(\ref{CC}) much as before but with these extra features and the same $J$. The Dirac operator then comes out as stated from
\begin{align*} \dirac e^\alpha&=\la\nabla_\CS e^\alpha=s^i\la S^\alpha{}_{i\beta}e^\beta=S^\alpha{}_{i\beta}C^{i\beta}{}_\gamma e^\gamma=O_{ij}({\imath\sigma^j\over 4}\mu_j+d_j)\sqrt{\lambda_j}O_{ik}{\sigma^k\over\sqrt{\lambda_k}}\\ &=\delta_{jk}({\imath\sigma^j\over 4}\mu_j+d_j)\sigma^k=d_j\sigma^j+{\imath\over 4}\sum_i
\mu_i=d_j\sigma^j+{\imath\over 4}{{\rm Tr}(g)\over \sqrt{\det(g)}}\end{align*}
and the Leibniz rule for $\nabla_\CS$. 

Finally, we use the same Hilbert space $\CH$ completing $\CS$ as before and the same inner product defined by $\int$. This works because the coefficients $g_{ij}$ are constants and hence $C^i$ is also a constant. So the proof of Proposition~\ref{dirround} is unchanged except that now we need $C^i$ to be hermitian not $\sigma^i$. But this holds as the $O_{ij}$ and $\sqrt{\lambda_j}$ are real under our assumption of Euclidean signature. Likewise ${\rm Tr}(g)/\sqrt{\det(g)}$ is real, so this term contributes correctly. As before, we need $d_i=0$ for self-adjointness of $\imath\dirac$. 

 \endproof

As before, this result holds up to a unitary conjugation in spinor space (because we had a unique solution for the preferred $C^i$ in the analysis above) and  phase parameter in $\CJ$. In component terms, our Dirac operator with $d_i=0$ is the formula (\ref{diracg}) as promised. Due to the geometric conditions (\ref{CSGam})--(\ref{CC}), we also know that the Lichnerowicz formula \cite[Prop. 8.45]{BegMa} applies relating $\dirac^2$ to the spinor Laplacian and curvature, i.e. this is a good Dirac operator from the quantum Riemannian geometry point of view as well as from the Connes one. 

\section{Concluding remarks and outlook}\label{secrem}

It is striking that making a geometric assumption -- notably (\ref{CSGam})-(\ref{CC}) for the Clifford structure, which is not needed for a spectral triple, i.e., the {\em conjunction} of quantum Riemannian geometry working from the bottom up and Connes' axioms working from the top down, led to a unique answer for up to unitary equivalence and an undetermined phase in $\CJ$. We emphasised the round metric case but the same applied for general $g_{ij}$  at least with Euclidean signature. That the KO dimension is $n=3$ mod 8 in the Euclidean signature case is in keeping with a `dimension jump' phenomenon in \cite{ArgMa} i.e., arising from the extra cotangent direction $\theta'$.  

The Lorentzian case with some of the $\lambda_i<0$ works the same way but should be studied further. From the geometric side, the derivation and formula for $\dirac$ in Proposition~\ref{propdirg} still applies but the last part of the proof concerning $\CJ$ and self-adjointness no longer apply in the form stated. For example, if we take the metric $g={\rm diag}(-1,1,1)$  with $\lambda_1=-1$ and $\lambda_2=\lambda_3=1$, the natural Clifford structure from (\ref{CO}) with $O=\id$ has $C^1=-\imath\sigma^1$, $C^2=\sigma^2$ and $C^3=\sigma^3$ while the QLC from  (\ref{QLC}) has 
\[ \Gamma^1{}_{jk}=\eps_{1jk},\quad \Gamma^2{}_{31}=-1,\quad \Gamma^2{}_{13}=-3,\quad \Gamma^3{}_{12}=3,\quad \Gamma^3{}_{21}=1\]
leading to a compatible spinor connection
\[ S_1=\imath({3\over 4}\sigma^1+ d_1),\quad S_2=-{1\over 4}\sigma^2+ d_2,\quad S_3=-{1\over 4}\sigma^3+ d_3\]
from (\ref{Sigen}) and $\dirac$ as in  Proposition~\ref{propdirg}. This time, however, $\CJ=q\sigma^3$ up to a phase $q$ solves (\ref{JJ}), (\ref{CJ}) with $\eps=\eps'=1$, while (\ref{SJ})  is solved by $d_2,d_3\in\R$ and $d_1\in \imath\R$. Thus,  the axioms of spectral triple that do not involve the Hilbert space hold (noting as before that we should flip the sign of $\eps'$ if we consider $\imath\dirac$ as the operator for the spectral triple). This time, however, neither $\dirac$ nor $\imath\dirac$ are self-adjoint with respect to our previous Hilbert space structure, which would need to be adapted. 

Our Dirac operator and spectral triple in the case of the `round metric' has a clean rotationally invariant form
\[  (\dirac\psi)_\alpha=\del_i\psi_\beta\sigma^{i\beta}{}_\alpha+{3\imath\over 4}\psi_\alpha\]
built from Pauli matrices and certain noncommutative vector fields $\del_i$ (in fact inner derivations) deforming Killing vectors for rotation about the three axes. It is in the same spirit as a previous ad-hoc but likewise rotationally invariant proposal \cite{BatMa} for a canonical Dirac operator on $\C_\lambda[\R^3]=U(su_2)$ as a fuzzy $\R^3$ with its 4D calculus with derivatives $\del^\alpha{}_\beta$. This is similar to ours considering the Pauli matrices as converting from vector to spinor indices, but was not a spectral triple nor obtained from a bimodule connection. Indeed, the 4D calculus used in \cite{BatMa} does not descend to $\C_\lambda[S^2]$.
 
Next, our conditions (\ref{JJ})--(\ref{CC}) all take place for constant coefficients; the algebra $A=\C_\lambda[S^2]$ does not enter. Hence, our results apply also when we set $\lambda_p=1/n$ and quotient out by the kernel of the $n$-dimensional representation. This quotient is isomorphic to $M_n(\C)$ via the $n$-dimensional representation and our differential and integral structures map over by 
\begin{equation}\label{reduction} x_i\mapsto {2\over n}J_i,\quad \del_i\mapsto -\imath[J_i,\ ],\quad \int\mapsto {1\over n}{\rm Tr}\end{equation}
 essentially as in \cite{Mad}. Here $J_i$ denotes the standard $n$-dimensional matrix representation of the $su_2$ relations $[J_i,J_j]=\imath\eps_{ijk}J_k$ and the Hilbert space norm is given by $\<a,b\>={1\over n}{\rm Tr}(a^\dagger b)$. In this way, our Dirac operator descends to a spectral triple on each matrix algebra reduced fuzzy sphere. Comparing with the subsequent literature, it is clear for the round metric that we  land essentially on the `full spectral triple' $\widetilde{\mathcal D}_{N}$ proposed in \cite{And} with $N=n-1$, which in turn was built on the Dirac operator of \cite{BalPad} for $j=1/2$ in their family. The Dirac operator here was motivated by $SU_2$-equivariance and constructed as the action of a `Dirac element' in $\mathcal D\in U(su_2)\tens U(su_2)$, whereas for us the rotational invariance is a derived property of the quantum geometry associated to the round quantum metric. \cite{And} also considers further aspects of Connes theory such as spectral distance, which it would be interesting to connect to the quantum Riemannian geometry. More precisely, \cite{And} studies this for their operator $D_N$ but with methods which can be adapted to $\widetilde{\mathcal D}_{N}$ for the reduced version of our $\imath\dirac$ for the round metric. For example, let  $\lambda_p=1/2$ which corresponds to reduced a fuzzy sphere with $n=2$, $N=1$, $j=1/2$. Then the difference is that $D_1$ in \cite{And} has  $\CH=\C^2\tens\C^2$ with $A=M_2(\C)$ acting in the vector representation on the first copy, whereas our $\imath\dirac$ has $\CH=M_2(\C)\tens\C^2$ with $A$ acting by left multiplication. 
 
\begin{proposition} $\imath\dirac$ with the round metric on the $\lambda_p=1/2$ reduced fuzzy sphere $A=M_2(\C)$ has the same spectral distance as $D_1$ in \cite{And}. In particular, for pure states viewed as coherent states associated to points in $S^2$ and separated by angle $\Theta\in[0,\pi]$, the spectral distance is $\sin(\Theta/2)$. \end{proposition}
\proof The spectral distance\cite{Con} we compute between two positive linear maps (states) $\omega,\omega':A\to\C$ is  
\[ d(\omega,\omega'):={\sup}_{a\in A\ :\ 
||[\imath\dirac, a]||\le 1}\left\{|\omega(a)-\omega'(a)|\right\}\]
which does not pay attention for the constant term of $\imath\dirac$. In this case $[\imath\dirac, a]=\imath\del_k a\tens \sigma^k{}^T={1\over 2}[\sigma^k,a]\tens\sigma^k{}^T$ (since we use column vectors for the second factor). Writing $a=a_i\sigma^i\in M_2(\C)$ since we can ignore any identity matrix component, we have $[\imath\dirac, a]=\imath a_i\eps_{ijk}\sigma^j\tens\sigma^k{}^T$ where the left copy acts by multiplication. We also have $[\imath\dirac, a]^\dagger=-[\imath\dirac,a^\dagger]$ with respect to the inner product on $\CH$ (which is via the trace on $M_2(\C)$ and the standard inner product on $\C^2$). Hence 
\begin{align*}[\imath\dirac, a]^\dagger&[\imath\dirac, a]= \bar a_l a_i \eps_{ijk} \eps_{lmn}\sigma^m\sigma^j\tens(\sigma^k\sigma^n) ^T\\
&= \bar a_l a_i\eps_{ijk} \eps_{lmn}(\delta_{mj}+\imath\eps_{mjp}\sigma^p)\tens(\delta_{kn}+\imath\eps_{kns}\sigma^s{}^T)
\\
&=2|a|^21\tens 1+\imath (\bar a\times a)\cdot(\sigma\tens 1-1\tens\sigma^T) + a\cdot\sigma\tens \bar a\cdot\sigma^T+ \bar a\cdot\sigma\tens a\cdot\sigma^T\end{align*}
using a 3-vector notation for $a=(a_i)$ and its complex conjugate, and $|a|^2=\bar a\cdot a$. We next chose the $1,\sigma^i$ basis of $M_2(\C)$ and write down matrices for the action of $\sigma^j$ by left multiplication. This converts the above to an $8\times 8$ matrix and we find eigenvalues $0,4|a|^2$. Hence $||[\imath\dirac,a]||=2|a|$ for the operator norm. This is the same $||[D_1,a]||$ in \cite{And} after which the calculation of $d(\omega,\omega')$ is the same as given there. As everything scales under scaling of $a$, it suffices for the supremum to take $|a|=1/2$. 

For completeness, we recall that a pure coherent state at $(\theta,\phi)$ for the azimuthal angle and rotation angle about the polar axis is given by  $\omega_{\theta,\vartheta}(a)=\< \theta,\varphi| a  |\theta,\varphi\>$ for 
 \[ |\theta,\varphi\>={1\over\sqrt{2}}\begin{pmatrix} e^{\imath\varphi}\sqrt{1-\cos\theta}\\ \sqrt{1+\cos\theta}\end{pmatrix}.\]
Such states  \cite{Are} are defined for all spins $j$ and such that $\omega_{\theta,\varphi}(x)=r n_{\theta,\varphi}$ where the unit vector $n_{\theta,\varphi}$ points to the location on $S^2\subset\R^3$. We use the conventions of \cite{BatMa}, where it is shown that they are characterised within the spin $j$ representation as states of minimum variance in position vector, and that they have expected radius $r=2\lambda_pj=1/2$ in our case. From this information and $a\cdot\sigma=2a\cdot x$ for $a=\alpha+\imath\beta$ regarded as a pair of real vectors $\alpha,\beta$ with $\sqrt{\alpha^2+\beta^2}=1/2$, it is easy to compute that
\begin{align*} |\omega_{\theta,\varphi}(a)-\omega_{\theta',\varphi'}(a)|&=|(\alpha+\imath\beta)\cdot(n_{\theta,\varphi}-n_{\theta',\varphi'})|\\
&= \sqrt{|\alpha\cdot(n_{\theta,\varphi}-n_{\theta',\varphi'})|^2 +|\beta\cdot(n_{\theta,\varphi}-n_{\theta',\varphi'})|^2}\\
&=|n_{\theta,\varphi}-n_{\theta',\varphi'}|\sqrt{\alpha^2\cos^2\theta_\alpha+\beta^2\cos^2\theta_\beta}\le {1\over 2}|n_{\theta,\varphi}-n_{\theta',\varphi'}| \end{align*}
for some angles $\theta_\alpha,\theta_\beta$ between $\alpha,\beta$ and the difference vector.  The bound is reached by $\theta_\alpha=\theta_\beta=0$ so $d(\omega_{\theta,\varphi},\omega_{\theta',\varphi'})=|n_{\theta,\varphi}-n_{\theta',\varphi'}|/2$ which comes out as $\sin(\Theta/2)$. The more sophisticated point of view in  \cite{And} argues the same formula also for mixed stated given by vectors $n$ of length possibly less than 1. 
\endproof

The spectral distances for $\imath\dirac$ for general $n$ as well as for the non-reduced fuzzy sphere will be looked at elsewhere. The proposal in \cite{Gro} also has some resemblance to our Dirac operator but does not construct a spectral triple, while \cite{Wat} were motivated by asking for $\dirac$ in a spectral triple to commute with a chirality operator (so as to have a symmetric spectrum). More broadly, \cite{Bar} defines any real finite spectral triple on a matrix algebra as a finite `fuzzy geometry' and identifies various types, including some with terms of the inner form $\sum_i\gamma^i\tens [L_i,\ ]$ for some $\gamma_i$ and some matrices $L_i$. Finally, for the reduction at $n=2$, our Dirac operator is a version for $\Omega^1$ 3-dimensional  of the Dirac operator on $M_2(\C)$ with its 2-dimensional calculus in \cite{BegMa:spe}, where $\dirac={1\over 2}(\sigma^1\tens [\sigma_1,\ ]- \sigma^2\tens [\sigma_2,\ ])$. 
 
There could nevertheless be other quantum-geometric spectral triples on the fuzzy sphere under different assumptions than the ones made here. First of all, we assumed a central basis $\{e^\alpha\}$ but one could also look more generally, for example relations of the form $[e^\alpha,x^i]=\lambda_p\sigma^i{}^\alpha{}_\beta e^\beta$ would be rotationally invariant and compatible  with the algebra relations in the sense of the Jacobi identity. This could lead to a larger moduli of geometric Dirac operators for the round metric. Thinking more geometrically, but still for the round metric, our spinor bundle $\CS$ was trivial with the spinor components elements of the `coordinate algebra' $A$ (which we understood as noncommutative spherical harmonics). This made sense since our cotangent bundle is also trivial, but it is not the only possibility. Indeed, we have already constructed a natural fuzzy monopole in \cite{LirMa} and we could look for Dirac operators based on this.  Here, the spinor bundle on the sphere would be of the form $\CS=\CS_+\oplus\CS_-$ where $\CS_+$ is the charge 1 monopole line bundle as found in \cite{LirMa} and $\CS_-$ is its dual. The Clifford action $\Omega^1\tens_A\CS\to \CS$ for the $q$-monopole and $q$-antimonopole case in \cite{BegMa:spe} was found from the holomorphic structure of its 2D calculus, which does not immediately apply here. This merits further study and could lead to a unified treatment including a fuzzy sphere Dirac operator with fuzzy monopole spherical harmonics, the  q-sphere one in \cite{BegMa:spe} and a $q$-fuzzy sphere one covering the 2-parameter spheres of Podl\`es in the unified formulation of \cite{Ma:fuz}, likely different from the isospectral q-deformation in \cite{Dab}. 

In terms of mathematical physics, quantum field theory on finite fuzzy spheres was studied in several works, e.g.  \cite{Gro1,Gro2}, but could be revisited for the unreduced case using the quantum geometric setting. Moreover, one could look at Dirac spinors on 4D quantum spacetime models such as the black hole model which a fuzzy sphere at each $r,t$ in \cite{ArgMa}, so as to explore particle creation with spinors and other issues. One could also revisit 2+1 quantum gravity models, i.e. on a fuzzy sphere at each time $t$ in the spirit of \cite{BatMa, FreMa,MaSch} for $\C_\lambda[\R^3]$ or in the context of spin networks\cite{Pen} and aspects of loop quantum gravity. It would also be interesting to see the content of the restriction to quantum-geometrically realised spectral triples in the context of Connes' approach to the standard model in \cite{Con0,ConMar} and related works such as \cite{Bar0, DabSit}.  These are some directions for further work.

\end{document}